Dmytro Taranovsky
November 5, 2003

# Determinacy Maximum

**Abstract:** We propose a new determinacy hypothesis for transfinite games, use the hypothesis to extend the perfect set theorem, prove relationships between various determinacy hypotheses, expose inconsistent versions of determinacy, and provide a philosophical justification for determinacy.

## Introduction

Many natural propositions that are undecidable in ZFC can be resolved by determinacy hypotheses. For example, projective determinacy provides a reasonably complete theory of second order arithmetic. This paper introduces a strong determinacy hypothesis, which we hope resolves more of the natural undecidable propositions.

First, we formalize the notion of the transfinite game and the strategies. The game is a perfect information game played by two players. The players move in turns, and at limit ordinals the first player moves. Each move consists of picking a member of the nonempty set of allowed moves, called the move set. A game on S means a game whose move set is S. A position is a sequence of moves at a point before the game ends. Once the sequence reaches a certain ordinal length, called the length of the game, the game ends. The play is the sequence of all moves selected by the players. If the sequence of moves belongs to the payoff set for the game (the payoff set is determined before the game starts), then the first player wins; otherwise, the second player wins.

A strategy for a player is a mapping that for every position for which it is the player's turn to move assigns a valid move. For a particular game, a strategy for one player is called winning if the player wins the game provided that he follows the strategy regardless of the choices of the other player. A game is said to be determined if one of the players has a winning strategy. A determinacy hypothesis asserts that each game in a particular class of games is determined.

Throughout this paper, we work in the standard set theory, ZFC. Real numbers are identified with subsets of ω. **R** is the set of all real numbers. A set is ordinal definable if there is an ordinal λ such that the set is definable in $V_\lambda$ by a first order formula with ordinal parameters. $OD_x$ is the class of sets that are ordinal definable using $x$ as a parameter.

## Determinacy Maximum

The following principle, which I shall call determinacy maximum, seems to be the strongest canonical determinacy hypothesis which is not (known to be) inconsistent:

**Determinacy Maximum:** Every game of length $\omega_1$ with ordinal definable payoff set and with the set of allowed moves being an ordinal is determined.

The obvious attempts of strengthening the principle are either inconsistent or trivial.

First, we recall the proof of the existence of undetermined sets.

**Theorem 1:** ZF proves that for every well-ordering X of the continuum with each element having less than the continuum predecessors, there is an undetermined game on {0, 1} of length $\omega$ that is definable (and ordinal definable) from X.

**Proof:** Choose a definable coding that for each strategy for the first and for the second player assigns a unique real number. Use X to define a well-ordering of the set of strategies with each strategy having less than $2^\omega$ predecessors. Construct a mapping $f$ from the set of real numbers into {0, 1, 2} as follows: Start with $f(x)=0$ for all $x$. Consider all of the strategies in the order specified by the well-ordering. At each point, the set of $x$ such that $f(x)\neq 0$ has cardinality less than the continuum. For each strategy, there are $2^\omega$ ways to play by the other player, each yielding a different sequence of moves. Thus, there is a sequence of moves $x$ with one player following the strategy such that $f(x)=0$. Pick the smallest such sequence according to the well-ordering X. If the strategy is for the first player, set $f(x)=2$; otherwise, set $f(x)=1$. Once all strategies have been examined, the construction of $f$ is completed, and set A = $\{x: f(x)=1\}$. A is definable from X. To see that A is not determined, consider a strategy for the first player. By construction, there is a sequence of moves x, with the first player following the strategy such that $f(x)=2$. Since x is not in A, the strategy is not winning. Similarly, for a strategy for the second player, there is a sequence of moves x, with the second player following the strategy such that $f(x)=1$. Since x is in A, the first player wins and the strategy is not winning.

**Theorem 2:** There is an undetermined game of length $\omega$ with ordinal definable payoff set and with allowed moves being arbitrary sets of real numbers.

**Proof:** We construct the game as follows. During the first move, the first player must pick an undetermined set X of real numbers, which is possible by theorem 1. Then, at every move, both players must play the empty set or **R** — the first player to violate the requirement loses. Assume that both satisfy the requirement. The sequence of moves, starting with the second move of the first player can be coded into a real number $r$. The first player wins if and only if X is undetermined and $r$ is in X. Since X is undetermined and since both players know the value of X before starting to pick the digits of $r$, the game is undetermined.

**Theorem 3:** There is an undetermined game on {0, 1} of length $\omega_1+\omega$ with ordinal definable payoff set.

**Proof:** The first $\omega_1$ moves of the first player code a set of real numbers X such that $x$ is in X iff there is a countable ordinal $\alpha$ such that the $\omega$ moves of the first player after the first $\omega*\alpha$ moves represent $x$. Clearly, the first player can code an arbitrary set of real numbers of cardinality not greater than $\omega_1$. Let $r$ be the play excluding the first $\omega_1$ moves.

Recall that for every X there is a definable (using X as a parameter) payoff set $G^*(X)$ such that the first player has a winning strategy if X has a perfect subset, the second player has a winning strategy if X is countable, and $G^*(X)$ is undetermined otherwise (Moschovakis 295).

We construct the game as follows:
If the continuum hypothesis is false, the first player wins iff X is uncountable and $r$ is in $G^*(X)$. Since X is neither countable nor contains a perfect subset, the game is undetermined. If the continuum hypothesis is true, the first player wins iff X is undetermined and $r$ is in X. Since X is undetermined, which is possible by theorem 1, and both players know its value beforehand, the game is undetermined.

**Conjecture:** If X is an uncountable sequence of distinct of real numbers, then there is an undetermined game on {0, 1} of length ω with the payoff set ordinal definable from X.
**Note:** A modification of the proof of theorem 3 easily yields: If X is an uncountable sequence of distinct of real numbers, then there is an undetermined game on {0, 1} of length ω with the payoff set ordinal definable from X and a real number.
**Update (Oct 29, 2016):** The conjecture is true. If X has length $\omega_1$ and a perfect subset, then a well-ordering of **R** with every element having countably many predecessors is constructible from X by picking the maximal perfect subset of X and mapping it to **R**.

Thus, using $2^\omega$ instead of $\omega_1$ in the definition of determinacy maximum would imply the continuum hypothesis. It seems that determinacy maximum is especially strong when it is combined with the continuum hypothesis, as then one can play an arbitrary set of real numbers, not merely a near countable one; one can even play an arbitrary mapping of the set of real numbers into ordinals.

If a game has length λ+n for an integer n, then once λ moves are made, the rest of the game becomes of finite length and hence determined, so determinacy for ordinal definable games of length λ+n is reduced to that of games of length λ, and without loss of generality one can limit consideration to games whose length is a limit ordinal. Theorem 3 shows that broad determinacy hypotheses are limited to games not longer than $\omega_1$.

It may be instructive to see how an attempt to disprove the combination of the determinacy maximum and the continuum hypothesis fails. We can force the first player to produce a well-ordering X of the real numbers, with each real number having countably many predecessors. We can also make such that the two players play an auxiliary game of length ω continuum many times and the first player wins only if each of the plays is in A(X), the undetermined set defined from X by the method specified in the rules. The second player cannot have a winning strategy since given a strategy, for each sub-game of length ω, the first player can choose a sequence of moves that defeats the strategy. It may seem that since A(X) is undetermined, the first player cannot guarantee a win either. However, X is not fully constructed until after all of the sub-games are played, and as each sub-game is played, the first player can fine-tune X so as to cause the previous sequences of moves to be in A(X). While one can impose conditions on X, there does not appear a way (without guaranteeing a winning strategy for the second player) to prevent the first player from choosing the order of elements of X so as to win in every sub-game. The above argument does suggest that determinacy maximum is a very strong statement about sets.

A modification of the proof of theorem 3 strengthens theorem 2.

**Theorem 4:** There is an undetermined game of length ω with ordinal definable payoff set and with allowed moves being sequences of zeros and ones of length $\omega_1$.

**Proof:** We construct the game as follows. The first move of the first player codes a set X of real numbers. For each move, starting with the second move of the first player, extract the first element, and collect the elements in order into a real number r. The first player wins iff X and r satisfy the requirement on the first player in the proof of theorem 3. The game is undetermined for the reasons stated in the proof of theorem 3.

Ordinal definability is a rather broad condition; it appears that theorems 1-4 can be sharpened by replacing it with a more narrow form of definability. Irregularities in properties of sets can make pathologies ordinal definable — for example, in building models one can choose for which successor cardinals the generalized continuum hypothesis holds, and the set of cardinals less than $\omega_2{\wedge}\omega$ that follow the generalized continuum hypothesis can code a well-ordering of real numbers. The broadness of definability in determinacy maximum implies, among other things, that the truth set of set theory is canonical in a certain way rather than random.

**Theorem 5:** Determinacy maximum implies that every game of length not greater than $\omega_1$ whose move set is an ordinal definable set of countable sequences of ordinals and whose payoff set is ordinal definable from a countable sequence of ordinals is determined.

**Proof:** Given a game whose moves are countable sequences of ordinals and whose payoff set is A, we reduce it to a game whose moves are ordinals and with the payoff set definable in terms of A. Partition $\omega_1$ into $\omega_1$ intervals, each of length $\omega$. Define the payoff set as follows: From interval number $\omega*\alpha+2n+1$, collect the moves of the first player into a countable sequence. From interval number $\omega*\alpha+2n+2$, collect the moves of the second player in order into a countable sequence. Collect the sequences in order into a sequence of length $\omega_1$. The first player wins if and only if the sequence belongs to A. Clearly, each winning strategy for the new game translates into a winning strategy for the original game, and determinacy of games with moves being countable sequences of ordinals is reduced to the corresponding determinacy of games with moves being ordinals.

Suppose that there is an undetermined game not longer than $\omega_1$ with moves being ordinals less than a particular ordinal and the payoff set ordinal definable from a countable sequence of ordinals. Then, there is an ordinal $\lambda$ such that there is an undetermined game not longer than $\omega_1$ with moves being ordinals less than $\lambda$, and with payoff set that is ordinal definable from a countable sequence of ordinals less than $\lambda$. We construct a game that contradicts the determinacy maximum as follows. The first $\omega$ moves of the first player must form a sequence *x* for which there is an undetermined game not longer than $\omega_1$ with moves being ordinals less than $\lambda$, and with the payoff set ordinal definable from *x*. Let X be the least undetermined payoff set (for games of length $\omega_1$) in the canonical well-ordering of OD$_\mathbf{x}$. The first player wins if and only if his first $\omega$ moves satisfy the requirement above and the sequence of moves after that is in X. Since X is undetermined and constructed before the sequence is played, the constructed ordinal definable game is undetermined.

The methods of the proof easily yield the following corollary:

**Corollary 6:** If all ordinal definable (that is the payoff set is ordinal definable) games of length $\omega+\alpha$ on {0, 1} are determined, then all games on {0, 1} ordinal definable from the reals (that is using real numbers as parameters) of length $\alpha$ are determined. If all ordinal definable games on {0, 1} of length $\omega*\alpha$ are determined, then all games on {0, 1} of length $\alpha$ whose payoff set is ordinal definable from the real numbers are determined.

**Note:** In corollary 6, "ordinal definable" can be replaced with a narrower notion of definability provided that the notion is sufficiently robust.

These theorems show that determinacy maximum maximizes the allowed length of the game, the set of allowed moves, and the class of payoff sets. The attempts to extend determinacy maximum either lead to an equivalent statement (theorem 5), or a contradiction (theorems 1, 2, 3, and 4). However, one can make ad hoc determinacy assumptions. For example (the idea for this is due to John Steel), one could allow the games to have an arbitrary length as long as the game terminates whenever the position is not definable from a countable sequence of ordinals.

Ordinary descriptive set theory deals with definable sets of sequences of integers of length ω; determinacy hypotheses for games of length ω allow full development of the theory. Determinacy maximum asserts determinacy for games of length $\omega_1$, so it should have important consequences for the theory on definable sets of sequences of length $\omega_1$. The consequences of determinacy maximum are especially strong when combined with the continuum hypothesis.

A tree is a set of sequences (under the inclusion ordering) that is closed under subsequences. For an infinite regular cardinal κ, a tree T of height κ is said to be perfect if every totally ordered subset of T is contained in a path (a totally ordered subset of T of order type κ), and every path has κ points at which the tree branches. A set of subsets of κ is a perfect subset if it is the set of paths through a perfect tree of height κ.

A somewhat weaker definition of a perfect subset — the stronger definition is the one used in this paper — would be as follows. Let $P(P(\omega_1))$ be the set of all sequences of length $\omega_1$ of zeros and ones. For an arbitrary countable sequence **x** of zeros and ones, let O(**x**) be the set of elements of $P(P(\omega_1))$ whose initial segment is **x**. The set of all O(**x**) forms a basis for the standard topology on $P(P(\omega_1))$. A nonempty subset of $P(P(\omega_1))$ would be called perfect if it is closed and contains no isolated points.

**Theorem 7:** Assume determinacy maximum. Every definable from a countable sequence of ordinals subset of $P(P(\omega_1))$ either has cardinality at most continuum or has a perfect subset.

**Proof:** The proof is similar to the proof from determinacy that every uncountable definable set of real numbers has a perfect subset. Let X a subset of $P(P(\omega_1))$ that is definable from a countable sequence of ordinals. The moves of the first player will be (possibly empty) finite or countable sequences of zeros and ones, and each move of the second player must be a zero or a one. The moves of the two players are collected in order into a sequence of zeros and ones of length $\omega_1$; the first player wins iff the sequence is in X. By theorem 5, the game is determined. If the first player has a winning strategy, then the set of points of $P(P(\omega_1))$ that can be reached with the player following the strategy is a perfect subset of X. It suffices to show that if the second player has a winning strategy S, then the cardinality of X is at most continuum. Let *x* be in X and the first player play as follows: Whenever possible, the player makes a move such that after the move of the second player according to S, the position corresponds to a subsequence of *x*. Since S is a winning strategy for the second player, at some position **r**, this will not be possible. We claim that *x* is uniquely determined by S and **r**. If the contrary were to hold, let α be the position of the least digit of *x* that is not determined by S and **r**. Because of the maximality of **r**, if after reaching **r**, the first player chooses to play the sequence such that the new position corresponds to the initial subsequence of *x* of length α, then the second player will play the negation of the digit number α of *x*. Because digits of *x* whose position is less than α are determined by **r** and S, digit number α is also determined by **r** and S, contrary to the assumption. Thus, there is a surjection of the set of positions into the set of elements of X. Since there are only continuum many positions, X has cardinality at most continuum.

An interesting question is whether every definable from $\omega_1$ ordinals subset of $P(P(\omega_1))$ either has cardinality at most $\omega_1$ or has a perfect subset. More generally,

**Question:** Is it consistent/true that for every infinite regular cardinal κ every definable from κ ordinals subset of P(P(κ)) either has cardinality at most κ or has a perfect subset?

**Note:** The consistency strength of the assertion is at least that of a proper class of inaccessible cardinals, and I conjecture that by generalizing Levy collapse one can force

from a ground model containing a proper class of inaccessible cardinals an affirmative answer. An affirmative answer in V would imply the generalized continuum hypothesis and the negation of the Kurepa Hypothesis for every infinite successor cardinal.

Thus, the determinacy maximum is useful, but the question remains whether it is consistent. And if so, whether it is true, which leads us to the next section.

# Philosophy of Determinacy Hypotheses

The dividing line between ordinary conjectures and proposed axioms is that ordinary conjectures are expected to be proved while proposed axioms tend to be provably unprovable. Moreover, if an axiom increases proof theoretical strength of the theory, as is usually the case with determinacy hypotheses, then the consistency of the axiom is unprovable from the consistency of the theory.

There are a priori and a posteriori reasons for accepting an axiom. An axiom may be seen true through careful examination of its statement and a contemplation of the mathematical universe. Alternatively, or in addition, the consequences of an axiom may form a coherent whole and resolve — in seemingly the right way — the multitude of undecidable propositions in a given field, which compels at least some acceptance of the axiom by the practitioners of the field irrespective of whether it has a priori justifications.

One way of finding axioms is by transferring our proven intuitions into the infinite. For finite games, determinacy is obviously true — assuming flawless play, either the first player can win, or he cannot, in which case the second player wins. It is also provable: Final positions are determined, a position is winnable if and only if there is a move that makes the position losable for the other player; by induction from last positions into previous ones, every position is winnable or losable.

The theorems above show that determinacy maximum, assuming that it is consistent, is the maximal canonical transfer of our intuition about games and strategies. The dividing line is not arbitrary: To define a strategy, one must first *define* the game. For finite games, *every* payoff set and position is defined. For infinite games, determinacy is thus asserted only for games that are definable in a certain way. To apply determinacy maximum, the payoff set must be definable. (If determinacy maximum is false, then the least undetermined game in a canonical well-ordering of ordinal definable sets is definable by a finite formula of reasonable size.) In addition, every position is definable using a countable set of ordinal parameters. At the end of the game, the only thing that matters is if the position belongs to the definable payoff set; the actual, possibly indefinable, sequence of moves is irrelevant. The boundary of countably many parameters is chosen because an uncountable set of parameters can be used to construct pathological objects — like decomposition of a sphere into finitely many pieces which can be reassembled after translation and rotation into two spheres each congruent to the original — objects that appear indefinable and unconstructible from less than continuum parameters. If the continuum hypothesis is true — as every natural uncountable set of real numbers has cardinality continuum — then $\omega_1$ parameters suffice for all such paradoxes. (An ordinal is considered as one parameter since the notion of an ordinal is simply the extension of the notion of an integer — one keeps counting past infinity — and is apparently incapable of coding arbitrary infinite sets of integers.) In the restriction of determinacy maximum to games whose moves must be 0 or 1 — and such restriction is still very strong — every position is defined by a real number — and every position can be recorded on a rule as a distance from mark zero.

The use of maximal canonical transfer of intuitions into the infinite is not new in set theory

— ZFC itself is the maximal canonical transfer of our basic intuitions on (finite) sets. The fact that the most obvious transfer of intuitions (every game is determined) is contradictory is also not new: Before ZFC was discovered, the full comprehension principle was found inconsistent by Russell. Human beings appear to lack direct access to infinite sets, so use of our indirect knowledge is necessary to obtain a good theory of them. The separation axioms are one piece of such knowledge; determinacy hypotheses are another.

Set theory is the theory of all sets, so it needs axioms to state the existence of arbitrary sets. However, there is no mathematical way to state that every set exists: a theory can only state closure under certain existence axioms. While such closure produces all finite sets, it cannot guarantee more than a countable number of infinite sets: In accordance with downward Löwenheim-Skolem theorem, every sound set theory has a countable transitive model.

However, one can make approximations to the claim that every set exists. For many concrete areas of mathematics, separation and replacement axioms suffice. Other areas need the axiom of choice. Once axioms state the existence of required sets, the theories can become practically complete, as is the case with arithmetic, with constructible sets, or in case of the axiom of determinacy, sets constructible from the reals. Without the required existence axioms, the key questions become undecidable, as is the case for projective sets in ZFC. While it is easy to claim existence of simple sets, there is no obvious way to claim existence of complicated sets, so standard axiom systems guarantee existence of the basic sets but are deficient on the complex ones.

Projective determinacy is so effective in second order arithmetic because it approximates the claim that arbitrary sets of integers exist. Rules for the games are often easily defined; false strategies usually have plenty of counter-strategies, but winning strategies can be difficult to find. The rules of chess have been known for centuries, thousands of strategies exist, but no one has shown a winning (or drawing) strategy. If we limit the universe of chess strategies to those that are constructible by a computer, then there is probably no winning strategy since every chess program has weaknesses that can be exploited to find a constructible counter-strategy. Yet, a correct strategy in all finite games exists; and claiming its existence is one way of asserting existence of integers whose value we have not constructed.

Similarly, if the set universe is limited to L, then many definable games have no winning strategies because the strategies are too complex to be constructible, and claiming their existence is an effective way of asserting existence of many complex sets of integers. The existence of the strategies is in fact provable from seemingly unrelated assertions of existence of complicated countable sets: Projective determinacy is equivalent to the claim that for every positive integer $n$ and real $r$, there is an iterable transitive model of ZFC containing $r$ and with $n$ Woodin cardinals. (The Woodin cardinals imply that the structure of the models is very complex.) Projective determinacy can be viewed as an assertion that there are many complicated countable sets. In fact, a wide variety of unrelated propositions of set theory contain a sufficient set existence component to imply projective determinacy. Determinacy maximum asserts that there are many complicated sets of ordinals, which code winning strategies for difficult games, like the one considered after theorem 3. Other determinacy hypotheses are also equivalent to existence of complicated sets: Determinacy for analytic sets is equivalent to the existence of a sharp (a certain complicated set of integers) for every real number. Determinacy hypotheses are existence axioms, like the separation axiom scheme or the axiom of choice.

Since, however, determinacy hypotheses are less obvious than the separation axioms or even the axiom of choice, one must examine their consequences before accepting them as axioms. Projective determinacy provides a very nice and practically complete theory of projective sets. Every projective set is measurable, has Baire property, and has perfect subset property. Pointclasses $\mathbf{\Pi}^1_{2n+1}$ and $\mathbf{\Sigma}^1_{2n}$ have the scale property, which is a very strong structural property and implies the properties of prewellordering and of uniformization. Proofs from projective determinacy are much more natural than proofs from V=L. Unlike semi-complete theories of projective sets obtained by restricting consideration to only certain not too complicated sets — the axiom of constructibility being the most prominent example — projective determinacy is unique in that it provides a canonical theory of projective sets without being restrictive. Other hypotheses that provide a nonrestrictive canonical theory of projective sets imply projective determinacy.

By claiming existence of complicated sets, determinacy hypotheses also have large proof theoretical strength. A theory cannot prove propositions of larger proof theoretical strength than itself. By limiting oneself to ZFC, one is cut off from such propositions — most of which are accessible through determinacy hypotheses.

Projective determinacy only claims determinacy for games where positions can be coded by integers and payoff sets definable in a simple way from real numbers. For the a priori and a posteriori reasons, the mathematical community should accept it, or a statement that implies it, as an axiom.

The consequences of determinacy maximum are yet to be explored. One hopes that they lead to a canonical theory of the larger fragments of the set theoretical universe, but first one must investigate the consistency of determinacy maximum.

REFERENCES

Moschovakis, Yiannis. *Descriptive Set Theory*. North-Holland Publishing Company: New York, 1980.
Woodin, W. Hugh. "The Continuum Hypothesis, Part I." *Notices of the AMS*. June/July 2001.

**Appendix:** Updates and New Results

**Update:** Feb 20, 2006
Determinacy maximum limited to games on zeroes and ones is consistent assuming that existence of a Woodin limit of Woodin cardinals and a measurable above them is consistent, with conjectured equiconsistency at Woodin limit of Woodin cardinals.

*A Strengthening of AD*
There is also a strong determinacy hypothesis for ZF + DC (inconsistent with AC):
Every game of length omega is weakly quasi-determined.

We define a weak quasi-strategy as follows. Before each move, the player to move chooses an arbitrary non-empty set, and receives an element of the set (just receives, not receives from the other player; this is not a part of the game and does not affect who wins). The moves may depend not only on the previous moves, but also on the elements received. To be winning, a strategy must work for all possibilities.

For example, if the first move is a countable ordinal α, the second player may choose the set of all reals of rank α, and will receive a particular real of rank α.

The hypothesis implies the axiom of determinacy. Under ZF + DC, at most one player has a winning weak quasi-strategy. In fact, if we add a generic well-ordering of the universe, winning weak quasi-strategies can be converted into winning strategies, leading (assuming enough cardinals will not be collapsed) to ZFC + every game on ordinals of length omega with the payoff set ordinal definable from a countable sequence of ordinals is determined.

**Update:** February 20, 2012
One consequence of determinacy maximum is that for every uncountable cardinal κ, there is a definable normal fine countably-complete ultrafilter on definable from countable sequence of ordinals sets of countable subsets of κ. In the relevant game, there are ω rounds, and at each round a player picks a countable subset of κ. The payoff considers only the union of all subsets played, and considered this way, the payoff belongs to the ultrafilter iff the first player has a winning strategy.

**Update:** March 7, 2012
Here is an extension of determinacy maximum. We do not know whether it is consistent.
**Extended Determinacy Maximum:**
Ordinal definable (in the sense below) two player perfect information games of arbitrary ordinal length with positions being sets definable from countable sets of ordinals are determined.
Here, ordinal definable means that there is ordinal definable:
- initial position
- move function: position, player number, move → new position
- limit function f (applied at all limit stages): f(p) → new position, where p is a sequence of positions corresponding to a valid play (starting at the initial position) and such that for every valid s and t of the same length f(s) = f(t) holds if for cofinally many α < length(s), s(α)=t(α).
- payoff set (set of positions won by the first player)

**Note:** One weakening is to use "ordinals" instead of "sets definable from countable sets of ordinals".

The strict requirements on the limit function are necessary to prevent indeterminacy for games of length $\omega_1 * \omega$. Otherwise, we could define a nondetermined game as follows: Force the first player to play a well-ordering of the reals on each $\omega_1$ round, and require that for all but finitely many rounds the well-ordering is the same. Use the well-ordering to get (in a determined way) a game on {0,1} of length omega with all positions undetermined, and make the two players play that game. The second player does not have a winning strategy, and given a strategy for the first player on the {0,1} game, have the second player defeat the strategy. This forces the first player to change the well-ordering infinitely many times and thus lose.

A consequence of the extended determinacy maximum is that for every regular uncountable κ, no ordinal definable canonical stationary subset of κ can be split in an ordinal definable way into disjoint stationary sets. For example, there are no disjoint ordinal definable stationary subsets s and t such that all members of s∪t have the same cofinality. This can be proved by considering games where players pick ordinals and the payoff depends only on the supremum of the ordinals picked. Determinacy maximum suffices for cofinalities ω and $\omega_1$.
**Note:** Part of the proof is to get a history-free winning strategy T', which given a winning strategy T is obtained as follows: Construct such T' from T by recursively mapping positions consistent with T into plays consistent with T, and doing this in a consistent way.

**Update:** December 14, 2016

Here are two consequences of determinacy maximum, even when limited to games on {0,1}:

- If φ is a formula with one free variable, and S and T are sufficiently fast-growing sequences of Turing degrees of length $\omega_1$, then φ(S)⇔φ(T).

- Assuming CH, all sufficiently high degrees of subsets of $\omega_1$ under constructible (or even recursive) reducibility are also indistinguishable by definable formulas.

For formal statements and proofs, see "Determinacy and Fast-growing Sequences of Turing Degrees", arXiv:1612.04494.